УДК 519.856

# Универсальный метод для задач стохастической композитной оптимизации

А.В. Гасников[1,2] gasnikov.av@mipt.ru

Ю.Е. Нестеров[3,4] yurii.nesterov@uclouvain.be

[1] Институт проблем передачи информации им. А.А. Харкевича Российской академии наук.
127051, Россия, г. Москва, Большой Каретный переулок, д.19 стр. 1
[2] Лаборатория структурных методов анализа данных в предсказательном моделировании (ПреМоЛаб), Факультет управления и прикладной математики Национального исследовательского Университета «Московский физико-технический институт».
141700, Россия, Московская область, г. Долгопрудный, Институтский переулок, д. 9
[3] Center for Operation Research and Econometrics Université Catholique de Louvain.
Voie du Roman Pays 34, L1.03.01 - B-1348 Louvain-la-Neuve (Belgium)
[4] Департамент больших данных и информационного поиска, факультет Компьютерных наук Национального исследовательского Университета «Высшая школа экономики».
125319, Россия, г. Москва, Кочновский проезд, 3

**Аннотация**

В работе предлагается быстрый градиентный метод для задач гладкой выпуклой оптимизации, требующий всего одну проекцию. Метод имеет наглядную геометрическую интерпретацию, поэтому получил название "метода подобных треугольников" (МПТ). В работе также предлагаются: композитный, адаптивный и универсальный вариант МПТ. Впервые (на базе МПТ) предлагается универсальный метод для сильно выпуклых задач (причем предложенный метод оказался непрерывным по параметру сильной выпуклости гладкой части функционала). Показывается, как универсальный вариант МПТ можно применять к задачам стохастической оптимизации.

**Ключевые слова:** быстрый градиентный метод, композитная оптимизация, универсальный метод, сильно выпуклый случай, стохастическая оптимизация, метод подобных треугольников.

## 1. Введение

В цикле работ [1 – 10] была предложена линейка быстрых градиентных методов для задач композитной оптимизации (выпуклый и сильно выпуклый случаи). Рассматривались приложения этих методов, в том числе к задачам стохастической оптимизации [3, 8 – 10]. Однако открытым оставался вопрос: возможно ли на базе этих методов предложить универсальный метод (настраивающийся на гладкость задачи) для задач стохастической оптимизации и(или) в сильно выпуклом случае? В данной работе получены ответы на эти вопросы. В основе предлагаемого подхода лежит оригинальный вариант классического быстрого градиентного метода (см. [1]): метод подобных треугольников (МПТ), изложенный в п. 2. Особенностью МПТ является необходимость выполнения всего одного "проектирования" на каждой итерации. Как следствие, предложенный метод оказался заметно проще многих своих аналогов. Это не только позволило упростить и привести в статье известные результаты для быстрых градиентных методов, но и продвинуться в решении отмеченных выше вопросов. В п. 3 МПТ распространяется на задачи сильной выпуклой композитной оптимизации. В п. 4 предлагается универсальный вариант МПТ для задач выпуклой и сильно выпуклой композитной оптимизации. В заключительном п. 5 универсальный МПТ из п. 4 переносится на задачи стохастической оптимизации.

## 2. Метод подобных треугольников для задач композитной оптимизации

Рассматривается задача выпуклой композитной оптимизации [2]
$$F(x) = f(x) + h(x) \to \min_{x \in Q}. \qquad (1)$$
Положим $R^2 = V(x_*, y^0)$, где прокс-расстояние (расстояние Брэгмана) определяется формулой (см., например, [1], главу 2 [3], [4, 11])
$$V(x,z) = d(x) - d(z) - \langle \nabla d(z), x - z \rangle;$$
прокс-функция $d(x) \geq 0$ ($d(y^0) = 0$) считается сильно выпуклой относительно выбранной нормы $\|\ \|$, с константой сильной выпуклости $\geq 1$; $x_*$ – решение задачи (1) (если решение не единственно, то выбирается то, которое доставляет минимум $V(x_*, y^0)$).

**Предположение 1.** *Для любых $x, y \in Q$ имеет место неравенство*
$$\|\nabla f(y) - \nabla f(x)\|_* \leq L \|y - x\|.$$

Опишем вариант быстрого градиентного метода для задачи (1) с одной "проекцией", который мы далее будем называть "метод подобных треугольников" (МПТ).

Положим

$$\varphi_0(x) = V(x, y^0) + \alpha_0 \left[ f(y^0) + \langle \nabla f(y^0), x - y^0 \rangle + h(x) \right],$$

$$\varphi_{k+1}(x) = \varphi_k(x) + \alpha_{k+1} \left[ f(y^{k+1}) + \langle \nabla f(y^{k+1}), x - y^{k+1} \rangle + h(x) \right], \quad (2)$$

$$A_k = \sum_{i=0}^{k} \alpha_i, \quad \alpha_0 = L^{-1}, \quad A_k = \alpha_k^2 L, \quad k = 0, 1, 2, \ldots, \quad x^0 = u^0 = \arg\min_{x \in Q} \varphi_0(x). \quad (3)$$

**Метод Подобных Треугольников**

$$\boxed{\begin{aligned} y^{k+1} &= \frac{\alpha_{k+1} u^k + A_k x^k}{A_{k+1}}, \\ u^{k+1} &= \arg\min_{x \in Q} \varphi_{k+1}(x), \\ x^{k+1} &= \frac{\alpha_{k+1} u^{k+1} + A_k x^k}{A_{k+1}}. \end{aligned}} \quad (4)$$

**Лемма 1 (см. [1]).** *Последовательность $\{\alpha_k\}$, определяемую формулой (3), можно задавать рекуррентно*

$$\alpha_{k+1} = \frac{1}{2L} + \sqrt{\frac{1}{4L^2} + \alpha_k^2}.$$

*При этом*

$$A_k \geq \frac{(k+1)^2}{4L}.$$

**Лемма 2.** *Пусть справедливо предположение 1. Тогда для любого $k = 0, 1, 2, \ldots$ имеет место неравенство*

$$A_k F(x^k) \leq \varphi_k^* = \min_{x \in Q} \varphi_k(x) = \varphi_k(u^k). \quad (5)$$

**Доказательство.** Проведем по индукции. При $k = 0$ формула (5) следует из того, что для любого $x \in Q$

$$F(x) = f(x) + h(x) \leq f(y^0) + \langle \nabla f(y^0), x - y^0 \rangle + \frac{L}{2} \|x - y^0\|^2 + h(x) \leq$$

$$\leq L V(x, y^0) + f(y^0) + \langle \nabla f(y^0), x - y^0 \rangle + h(x) = \frac{1}{A_0} \varphi_0(x).$$

Последнее неравенство следует из того, что для любых $x, z \in Q$

$$V(x,z) \geq \frac{1}{2}\|x-z\|^2.$$

Это неравенство, в свою очередь, следует из $\geq 1$-сильной выпуклости $d(x)$ в норме $\|\ \|$.

Итак, пусть формула (5) установлена для $k$, покажем, что тогда она будет справедлива и для $k+1$. По определению (2)

$$\varphi_{k+1}^* = \min_{x \in Q} \varphi_{k+1}(x) = \varphi_{k+1}(u^{k+1}) =$$

$$= \varphi_k(u^{k+1}) + \alpha_{k+1}\left[ f(y^{k+1}) + \langle \nabla f(y^{k+1}), u^{k+1} - y^{k+1}\rangle + h(u^{k+1})\right]. \tag{6}$$

Поскольку по предположению индукции $A_k F(x^k) \leq \varphi_k(u^k)$, и $\varphi_{k+1}(x)$ – сильно выпуклая в $\|\ \|$-норме функция с константой $\geq 1$ (это следует из аналогичного свойства функции $V(x, y^0)$, что, в свою очередь, следует из аналогичного свойства функции $d(x)$), то

$$\varphi_k(u^{k+1}) \geq \varphi_k(u^k) + \frac{1}{2}\|u^{k+1} - u^k\|^2 \geq A_k \cdot \left(f(x^k) + h(x^k)\right) + \frac{1}{2}\|u^{k+1} - u^k\|^2.$$

Из выпуклости $f(x)$ отсюда имеем

$$\varphi_k(u^{k+1}) \geq A_k f(y^{k+1}) + \langle \nabla f(y^{k+1}), A_k \cdot (x^k - y^{k+1})\rangle + A_k h(x^k) + \frac{1}{2}\|u^{k+1} - u^k\|^2. \tag{7}$$

Подставляя (7) в (6), получим

$$\varphi_{k+1}^* \geq A_{k+1} \cdot \underbrace{\left(\frac{A_k}{A_{k+1}} h(x^k) + \frac{\alpha_{k+1}}{A_{k+1}} h(u^{k+1})\right)}_{\geq h(x^{k+1})} + A_{k+1} f(y^k) +$$

$$+ \langle \nabla f(y^{k+1}), \underbrace{\alpha_{k+1} \cdot (u^{k+1} - y^{k+1}) + A_k \cdot (x^k - y^{k+1})}_{= A_{k+1} \cdot (x^{k+1} - y^{k+1})}\rangle + \underbrace{\frac{1}{2}\|u^{k+1} - u^k\|^2}_{= \frac{A_{k+1}^2}{2\alpha_{k+1}^2}\|x^{k+1} - y^{k+1}\|^2}. \tag{8}$$

Исходя из выпуклости функции $h(x)$ и описания МПТ (4), формулу (8) можно переписать следующим образом

$$\varphi_{k+1}^* \geq A_{k+1}\left[ f(y^{k+1}) + \langle \nabla f(y^{k+1}), x^{k+1} - y^{k+1}\rangle + \frac{A_{k+1}}{2\alpha_{k+1}^2}\|x^{k+1} - y^{k+1}\|^2 + h(x^{k+1})\right]. \tag{9}$$

Из предположения 1 следует, что если $A_{k+1}/\alpha_{k+1}^2 \geq L$, то

$$f\left(y^{k+1}\right)+\left\langle \nabla f\left(y^{k+1}\right), x^{k+1}-y^{k+1}\right\rangle +\frac{A_{k+1}}{2\alpha_{k+1}^{2}}\left\Vert x^{k+1}-y^{k+1}\right\Vert ^{2}\geq f\left(x^{k+1}\right). \qquad (10)$$

Согласно (3) $A_{k+1}/\alpha_{k+1}^{2}=L$, поэтому формула (10) имеет место. С помощью формулы (10) формулу (9) можно переписать следующим образом

$$\varphi_{k+1}^{*}\geq A_{k+1}\left[f\left(x^{k+1}\right)+h\left(x^{k+1}\right)\right]=A_{k+1}F\left(x^{k+1}\right).$$

Таким образом, шаг индукции установлен. Следовательно, лемма 2 доказана. ∎

Из лемм 1, 2 получаем следующий результат, означающий, что МПТ сходится как обычный быстрый градиентный метод, см., например, [1, 2–4, 11], т.е. МПТ сходится оптимальным образом для рассматриваемого класса задач [12].

**Теорема 1.** *Пусть справедливо предположение 1. Тогда МПТ (2) – (4) для задачи (1) сходится согласно оценке*

$$F\left(x^{N}\right)-\min_{x\in Q}F\left(x\right)\leq\frac{4LR^{2}}{\left(N+1\right)^{2}}. \qquad (11)$$

**Доказательство.** Из леммы 2 следует, что (в третьем неравенстве используется выпуклость функции $f(x)$)

$$A_{N}F\left(x^{N}\right)\leq\min_{x\in Q}\left\{ V\left(x,y^{0}\right)+\sum_{k=0}^{N}\alpha_{k}\left[f\left(y^{k}\right)+\left\langle \nabla f\left(y^{k}\right),x-y^{k}\right\rangle +h(x)\right]\right\} \leq$$

$$\leq V\left(x_{*},y^{0}\right)+\sum_{k=0}^{N}\alpha_{k}\underbrace{\left[f\left(y^{k}\right)+\left\langle \nabla f\left(y^{k}\right),x_{*}-y^{k}\right\rangle +h\left(x_{*}\right)\right]}_{\leq f(x_{*})+h(x_{*})}\leq$$

$$\leq V\left(x_{*},y^{0}\right)+\sum_{k=0}^{N}\alpha_{k}F\left(x_{*}\right)=R^{2}+A_{N}F\left(x_{*}\right). \qquad (12)$$

Заметим, что из второго неравенства следует, что если решение задачи (1) $x_*$ не единственно, то можно выбирать то, которое доставляет минимум $V\left(x_*, y^0\right)$. Именно таким образом возникает $R^2$ в оценке (12). Для доказательства теоремы осталось подставить нижнюю оценку на $A_N$ из леммы 1 в формулу (12). ∎

**Замечание 1.** В действительности в формуле (12) содержится более сильный результат, чем в формуле (11). А именно, формула (12) еще означает, что МПТ – прямо-двойственный метод (см., например, [13, 14]). Мы не будем здесь подробно на этом останавливаться, отметим лишь, что это свойство позволяет получать эффективные критерии останова для МПТ. Критерий останова позволяет не делать предписанного формулой (11) числа итераций и останавливаться раньше (по достижению желаемой точности). Это замечание можно распространить и на все последующее изложение.

МПТ получил такое название, поскольку на базе трех точек $u^k$, $u^{k+1}$, $x^k$ можно построить треугольник, а точки $y^{k+1}$ и $x^{k+1}$ лежат на сторонах этого треугольника ($y^{k+1}$ лежит на стороне $u^k x^k$, а $x^{k+1}$ на стороне $u^{k+1} x^k$), причем прямая, проходящая через точки $y^{k+1}$ и $x^{k+1}$, параллельна прямой, проходящей через точки $u^k$ и $u^{k+1}$.

Из описанной геометрической интерпретации получаются следующие результаты, распространимые и на все последующее изложение (с некоторыми оговорками в п. 5 [9])

**Следствие 1.** *Для любого $k = 0,1,2,...$ имеют место неравенства*

$$\left\| u^k - x_* \right\|^2 \leq 2V\left(x_*, y^0\right),$$

$$\max\left\{ \left\| x^k - x_* \right\|^2, \left\| y^k - x_* \right\|^2 \right\} \leq 4V\left(x_*, y^0\right) + 2\left\| x^0 - y^0 \right\|^2 = \tilde{R}^2.$$

**Доказательство.** Второе неравенство следует из первого, описанной выше геометрической интерпретации и неравенства $\left\| a+b \right\|^2 \leq 2\left\| a \right\|^2 + 2\left\| b \right\|^2$.

Докажем первое неравенство. Для этого воспользуемся леммой 1 и тем, что $\varphi_k(x)$ – сильно выпуклая в $\|\;\|$-норме функция с константой $\geq 1$. Для любого $x \in Q$

$$A_k F\left(x^k\right) + \frac{1}{2}\left\| x - u^k \right\|^2 \leq \varphi_k^* + \frac{1}{2}\left\| x - u^k \right\|^2 \leq \varphi_k(x) \leq$$

$$\leq \underbrace{\sum_{i=0}^{k} \alpha_i \left[ f\left(y^i\right) + \left\langle \nabla f\left(y^i\right), x - y^i \right\rangle + h(x) \right]}_{\leq A_k F(x)} + V\left(x, y^0\right) \leq A_k F(x) + V\left(x, y^0\right).$$

Выбирая $x = x_*$ и используя то, что $F\left(x^k\right) \geq F\left(x_*\right)$, получим первое неравенство следствия 1. ∎

Ранее такого типа результаты (как следствие 1) устанавливались только для евклидовой прокс-структуры в некомпозитном случае [14]. Следствие 1 играет важную роль в случае неограниченных множеств $Q$, на которых нельзя равномерно ограничить константу $L$, поскольку гарантирует, что как бы "плохо" себя не вела функция $F(x)$ вне шара конечного радиуса $\tilde{R}$ (зависящего от качества начального приближения) с центром в решении $x_*$, это никак не скажется на скорости сходимости метода, поскольку итерационный процесс никогда не выйдет за пределы этого шара.

**Следствие 2.** *Пусть $h(x) \equiv 0$ (т.е. $F(x) = f(x)$) и $\nabla f(x_*) = 0$. Тогда*

$$\max\left\{F(x^N), F(y^N), F(z^N)\right\} - \min_{x \in Q} F(x) \le \frac{L\tilde{R}^2}{N^2}.$$

### 3. Метод подобных треугольников для сильно выпуклых задач композитной оптимизации

В данном пункте будем считать, что $f(x)$ в задаче (1) обладает следующим свойством.

**Предположение 2.** $f(x)$ – $\mu$-*сильно выпуклая функция в норме* $\|\ \|$, *т.е. для любых* $x, y \in Q$ *имеет место неравенство*

$$f(y) + \langle \nabla f(y), x - y \rangle + \frac{\mu}{2}\|x - y\|^2 \le f(x). \quad (13)$$

Введем (в евклидовом случае $\tilde{\omega}_n = 1$)

$$\tilde{\omega}_n = \sup_{x, y \in Q} \frac{2V(x, y)}{\|y - x\|^2}.$$

**Замечание 2.** Из работы [15] следует, что в сильно выпуклом случае (когда сильно выпукла гладкая часть функционала $f(x)$, как в предположении 2) естественно выбирать именно евклидову норму и прокс-структуру, т.е. в большинстве случаев можно считать $\tilde{\omega}_n = 1$. Поясним это примером из [15]. Число обусловленности (отношение константы Липшица градиента $L_1(f)$ к константе сильной выпуклости $\mu_1(f)$), например, для квадратичных функций

$$f(x) = \frac{1}{2}x^T A x - b^T x,\ x \in \mathbb{R}^n,$$

посчитанное, скажем, в 1-норме не может быть меньше $n$. В то время как число обусловленности, посчитанное в евклидовой норме, может при этом равняться 1. Действительно, пусть $\xi = (\xi_1, \ldots, \xi_n)$, где $\xi_k$ – независимые одинаково распределенные случайные величины $P(\xi_k = 1/n) = P(\xi_k = -1/n) = 1/2$. Тогда, учитывая, что $\|\xi\|_1 \equiv 1$,

$$\mu_1(f) \le E_\xi\left[\xi^T A \xi\right] = \frac{1}{n^2}\operatorname{tr}(A) \le \frac{1}{n}\max_{i,j=1,\ldots,n}|A_{ij}| = \frac{1}{n}L_1(f).$$

Кроме того, из замечания 2 следует также и то, что если $\|\ \| = \|\ \|_1$, то $\tilde{\omega}_n \ge n$. Это обстоятельство также говорит в пользу выбора евклидовой прокс-структуры. Тем не менее, для общности далее мы будем допускать, что используется прокс-структура отличная от евклидовой.

Перепишем формулы (2), (3) следующим образом ($\tilde{\mu} = \mu/\tilde{\omega}_n$)

$$\varphi_0(x) = V(x, y^0) + \alpha_0 \left[ f(y^0) + \langle \nabla f(y^0), x - y^0 \rangle + \tilde{\mu} V(x, y^0) + h(x) \right],$$

$$\varphi_{k+1}(x) = \varphi_k(x) + \alpha_{k+1} \left[ f(y^{k+1}) + \langle \nabla f(y^{k+1}), x - y^{k+1} \rangle + \tilde{\mu} V(x, y^k) + h(x) \right], \quad (14)$$

$$A_k = \sum_{i=0}^{k} \alpha_i, \quad \alpha_0 = L^{-1}, \quad A_{k+1} \cdot (1 + A_k \tilde{\mu}) = \alpha_{k+1}^2 L, \quad k = 0, 1, 2, \ldots \quad (15)$$

Сам метод по-прежнему будет иметь вид (4) с $x^0 = u^0 = \arg\min_{x \in Q} \varphi_0(x)$.

**Лемма 3 (см. [7]).** *Последовательность $\{\alpha_k\}$, определяемую формулой (15), можно задавать рекуррентно*

$$\alpha_{k+1} = \frac{1 + A_k \tilde{\mu}}{2L} + \sqrt{\frac{1 + A_k \tilde{\mu}}{4L^2} + \frac{A_k \cdot (1 + A_k \tilde{\mu})}{L}}, \quad A_{k+1} = A_k + \alpha_{k+1}. \quad (16)$$

*При этом*

$$A_k \geq \frac{1}{L}\left(1 + \frac{1}{2}\sqrt{\frac{\tilde{\mu}}{L}}\right)^{2k} \geq \exp\left(\frac{k}{2}\sqrt{\frac{\tilde{\mu}}{L}}\right).$$

**Лемма 4.** *Пусть справедливы предположения 1, 2. Тогда для любого $k = 0, 1, 2, \ldots$ имеет место неравенство*

$$A_k F(x^k) \leq \varphi_k^* = \min_{x \in Q} \varphi_k(x) = \varphi_k(u^k).$$

**Доказательство.** Доказательство аналогично доказательству леммы 2. В основе лежит неравенство

$$V(x, y^k) \geq \frac{1}{2}\|x - y^k\|^2,$$

с помощью которого ключевая формула (8) перепишется следующим образом

$$\varphi_{k+1}^* \geq A_{k+1} \cdot \left( \frac{A_k}{A_{k+1}} h(x^k) + \frac{\alpha_{k+1}}{A_{k+1}} h(u^{k+1}) \right) + A_{k+1} f(y^k) +$$

$$+ \langle \nabla f(y^{k+1}), \alpha_{k+1} \cdot (u^{k+1} - y^{k+1}) + A_k \cdot (x^k - y^{k+1}) \rangle + \frac{(1 + A_k \tilde{\mu})}{2}\|u^{k+1} - u^k\|^2.$$

Отличие от формулы (8) в следующем

$$\frac{1}{2}\|u^{k+1} - u^k\|^2 \to \frac{(1 + A_k \tilde{\mu})}{2}\|u^{k+1} - u^k\|^2.$$

Рассуждая дальше точно также как при доказательстве леммы 2, получим

$$\varphi_{k+1}^* \geq A_{k+1} \left[ f\left(y^{k+1}\right) + \left\langle \nabla f\left(y^{k+1}\right), x^{k+1} - y^{k+1} \right\rangle + \frac{A_{k+1} \cdot (1 + A_k \tilde{\mu})}{2\alpha_{k+1}^2} \left\| x^{k+1} - y^{k+1} \right\|^2 + h\left(x^{k+1}\right) \right]. \quad (17)$$

Из предположения 1 следует, что если $A_{k+1} \cdot (1 + A_k \tilde{\mu}) / \alpha_{k+1}^2 \geq L$, то

$$f\left(y^{k+1}\right) + \left\langle \nabla f\left(y^{k+1}\right), x^{k+1} - y^{k+1} \right\rangle + \frac{A_{k+1} \cdot (1 + A_k \tilde{\mu})}{2\alpha_{k+1}^2} \left\| x^{k+1} - y^{k+1} \right\|^2 \geq f\left(x^{k+1}\right). \quad (18)$$

Согласно (15) $A_{k+1} \cdot (1 + A_k \tilde{\mu}) / \alpha_{k+1}^2 = L$, поэтому формула (18) имеет место. С помощью формулы (18) формулу (17) можно переписать следующим образом

$$\varphi_{k+1}^* \geq A_{k+1} \left[ f\left(x^{k+1}\right) + h\left(x^{k+1}\right) \right] = A_{k+1} F\left(x^{k+1}\right). \blacksquare$$

Из лемм 3, 4 получаем следующий результат, означающий, что МПТ в сильно выпуклом случае сходится как обычный быстрый градиентный метод (с двумя проекциями), т.е. МПТ сходится оптимальным образом для рассматриваемого класса задач.

**Теорема 2.** *Пусть справедливы предположения 1, 2. Тогда МПТ (14), (15), (4) для задачи (1) сходится согласно оценке*

$$F\left(x^N\right) - \min_{x \in Q} F(x) \leq LR^2 \exp\left(-\frac{N}{2}\sqrt{\frac{\tilde{\mu}}{L}}\right). \quad (19)$$

**Доказательство.** Из леммы 4 следует, что (в третьем неравенстве используется то, что $\tilde{\mu} = \mu / \tilde{\omega}_n$ и сильная выпуклость функции $f(x)$ – см. формулу (13) предположения 2)

$$A_N F\left(x^N\right) \leq \min_{x \in Q} \left\{ V\left(x, y^0\right) + \sum_{k=0}^N \alpha_k \left[ f\left(y^k\right) + \left\langle \nabla f\left(y^k\right), x - y^k \right\rangle + \tilde{\mu} V\left(x, y^k\right) + h(x) \right] \right\} \leq$$

$$\leq V\left(x_*, y^0\right) + \sum_{k=0}^N \alpha_k \underbrace{\left[ f\left(y^k\right) + \left\langle \nabla f\left(y^k\right), x_* - y^k \right\rangle + \frac{\mu}{2}\left\| x_* - y^k \right\|^2 + h(x_*) \right]}_{\leq f(x_*) + h(x_*)} \leq$$

$$\leq V\left(x_*, y^0\right) + \sum_{k=0}^N \alpha_k F(x_*) = R^2 + A_N F(x_*). \quad (20)$$

Для того чтобы получить оценку (19) осталось подставить нижнюю оценку на $A_N$ из леммы 3 в формулу (20). ∎

В действительности, выше установлено более сильное утверждение.

**Теорема 3.** *Пусть справедливы предположения 1, 2. Тогда МПТ (14), (16), (4) для задачи (1) сходится согласно оценке*

$$F(x^N) - \min_{x \in Q} F(x) \le \min\left\{\frac{4LR^2}{(N+1)^2}, LR^2 \exp\left(-\frac{N}{2}\sqrt{\frac{\mu}{L\tilde{\omega}_n}}\right)\right\}. \qquad (21)$$

Теорема 3 означает, что МПТ (14), (16), (4) непрерывен по параметру $\mu$. К сожалению, при этом в (16) явно входит этот параметр $\mu$. Если значение этого параметра неизвестно, то с помощью рестартов (см., например, [16, 17]) можно получить оценку (21) увеличив константы не более чем в 4 раза, т.е. число вычислений градиента $\nabla f(x)$ (обычно именно это является самым затратным в шаге), необходимых для достижения заданной точности, увеличится не более чем в 4 раза.

Между сильно выпуклым и просто выпуклым случаями имеется глубокая связь, позволяющая, например, получить оценку (19) с помощью оценки (11) и наоборот. Другими словами, имея эффективные алгоритмы решения выпуклых / сильно выпуклых задач, можно предложить на их базе алгоритмы решения сильно выпуклых / выпуклых задач. Покажем это, следуя, например, [10] (приводимые далее конструкции давно и хорошо известны, мы здесь их приводим для полноты изложения).

Введем семейство $\mu$-сильно выпуклых в норме $\|\ \|$ задач ($\mu > 0$)

$$F^\mu(x) = F(x) + \mu V(x, y^0) \to \min_{x \in Q}. \qquad (22)$$

**Теорема 4.** *Пусть*

$$\mu \le \frac{\varepsilon}{2V(x_*, y^0)} = \frac{\varepsilon}{2R^2}, \qquad (23)$$

*и удалось найти $\varepsilon/2$-решение задачи (22), т.е. нашелся такой $x^N \in Q$, что*

$$F^\mu(x^N) - F_*^\mu \le \varepsilon/2.$$

*Тогда*

$$F(x^N) - \min_{x \in Q} F(x) = F(x^N) - F_* \le \varepsilon.$$

**Доказательство.** Действительно,

$$F(x^N) - F_* \le F^\mu(x^N) - F_* \le F^\mu(x^N) - F_*^\mu + \varepsilon/2 \le \varepsilon.$$

Здесь использовалось определение $F_*^\mu$ и формула (23)

$$F_*^\mu = \min_{x \in Q}\{F(x) + \gamma V(x, y^0)\} \le F(x_*) + \gamma V(x_*, y^0) \le F_* + \varepsilon/2. \blacksquare$$

Приведем в некотором смысле обратную конструкцию.

**Теорема 5.** *Пусть функция $F(x)$ — $\mu$-сильно выпуклая в норме $\|\ \|$. Пусть точка $x^{\bar{N}}(y^0)$ выдается МПТ (2) – (4), стартующим с точки $y^0$, после*

$$\bar{N} = \sqrt{\frac{8L\omega_n}{\mu}} \qquad (24)$$

*итераций, где*

$$\omega_n = \sup_{x \in Q} \frac{2V(x, y^0)}{\|x - y^0\|^2}.$$

*Положим*

$$\left[ x^{\bar{N}}(y^0) \right]^1 = x^{\bar{N}}(y^0)$$

*и определим по индукции*

$$\left[ x^{\bar{N}}(y^0) \right]^{k+1} = x^{\bar{N}}\left( \left[ x^{\bar{N}}(y^0) \right]^k \right), \ k = 1, 2, \ldots .$$

*При этом на $k+1$ перезапуске (рестарте) также корректируется прокс-функция (считаем, что так определенная функция корректно определена на $Q$ с сохранением свойства сильной выпуклости)*

$$d^{k+1}(x) = d\left( x - \left[ x^{\bar{N}}(y^0) \right]^k + y^0 \right) \geq 0,$$

*чтобы*

$$d^{k+1}\left( \left[ x^{\bar{N}}(y^0) \right]^k \right) = 0.$$

*Тогда*

$$F\left( \left[ x^{\bar{N}}(y^0) \right]^k \right) - F_* \leq \frac{\mu \|y^0 - x_*\|^2}{2^{k+1}}. \tag{25}$$

**Доказательство.** МПТ (2) – (4) согласно теореме 1 (см. формулу (11)) после $\bar{N}$ итераций выдает такой $x^{\bar{N}}$, что

$$\frac{\mu}{2} \|x^{\bar{N}} - x_*\|^2 \leq F(x^{\bar{N}}) - F_* \leq \frac{4LV(x_*, y^0)}{\bar{N}^2}.$$

Отсюда имеем

$$\|x^{\bar{N}} - x_*\|^2 \leq \frac{8LV(x_*, y^0)}{\mu \bar{N}^2} \leq \frac{1}{2} \|y^0 - x_*\|^2 \frac{8L\omega_n}{\mu \bar{N}^2}.$$

Поскольку

$$\bar{N} = \sqrt{\frac{8L}{\mu} \omega_n},$$

то

$$\|x^{\bar{N}} - x_*\|^2 \leq \frac{1}{2} \|y^0 - x_*\|^2.$$

Повторяя эти рассуждения, по индукции получим

$$F\left( \left[ x^{\bar{N}}(y_0) \right]^k \right) - F_* \leq \left( \frac{1}{2} \right)^k \|y^0 - x_*\|^2 \frac{4L\omega_n}{\bar{N}^2} = \frac{\mu \|y^0 - x_*\|^2}{2^{k+1}}. \ \blacksquare$$

**Замечание 3.** Оценка (24), (25) в итоге получается похожей на оценку (19). Однако в подходе, описанном в теореме 5, всегда можно добиться, чтобы $\omega_n = \mathrm{O}(\ln n)$, в том числе и в случае $\|\ \| = \|\ \|_1$ [15] (следует сравнить в этом случае с оценкой $\tilde{\omega}_n \geq n$, приведенной выше). Кроме того, предложенная в теореме 5 конструкция позволяет рассматривать более общий класс сильно выпуклых задач, в которых $f(x)$ уже не обязательно – сильно выпуклая функция (см. предположение 2), что позволяет избавиться в оценках числа обусловленности $L/\mu$ от ограничения, описанного в замечании 2. В ряде приложений эти степени свободы оказываются чрезвычайно важными [10]. Однако, как показывают численные эксперименты (проводимые в случае евклидовой прокс-структуры и с априорно известной точной оценкой параметра $\mu$), подход, описанный в теореме 5, может проигрывать в скорости МПТ (14), (16), (4) один-два порядка, т.е. для достижения той же точности подход из теоремы 5 может потребовать до 100 раз больше арифметических операций. Из оценок это никак не следует. Но если метод из теоремы 5 работает по полученным верхним оценкам, то МПТ (14), (16), (4) (с критерием останова связанным с контролем малости нормы градиентного отображения [3, 18]) из-за отсутствия рестартов (т.е. необходимости делать предписанное число итераций) может остановиться раньше, что и происходит на практике. К тому же МПТ (14), (16), (4) еще и непрерывен по параметру $\mu$.

Из написанного выше, кажется, что в процедуре рестартов (теорема 5) можно использовать вместо предписанного числа итераций $\bar{N}$ на каждом рестарте какой-нибудь критерий останова. В частности, дожидаться, когда норма (или квадрат нормы) градиента (а в общем случае, когда минимум достигается не в точке экстремума, – норма градиентного отображения) уменьшиться вдвое. С таким критерием останова нет необходимости делать предписанного числа итераций на каждом рестарте. Однако пока не известен способ рассуждений, который позволял бы показать, что такая процедура сохраняет при перенесении оптимальность оценок (метод, работающий оптимально не в сильно выпуклом случае, порождает оптимальный метод и в сильно выпуклом случае). Впрочем, имеются различные эффективные на практике способы более раннего выхода с каждого рестарта (см., например, [19]), позволяющие (в случае задач безусловной оптимизации с евклидовой прокс-структурой и отсутствием композитного члена) ускорить описанную выше конструкцию на порядок.

### 4. Универсальный метод подобных треугольников

В ряде приложений значение константы $L$, необходимой МПТ для работы (см. формулу (16)), не известно. Однако, как следует из формул (10), (18), знание константы $L$ необязательно, если разрешается на одной итерации запрашивать значение функции в нескольких точках. Опишем соответствующий адаптивный вариант МПТ (14), (16), (4) (АМПТ) (см., например, [2]).

Положим $A_0 = \alpha_0 = 1/L_0^0$ – параметр метода (считаем здесь и везде в дальнейшем $L_0^0 \le L$, иначе во всех приводимых далее оценках следует полагать $L := \max\{L_0^0, L\}$),

$$k = 0, \; j_0 = 0; \; x^0 = u^0 = \arg\min_{x \in Q} \varphi_0(x).$$

До тех пор пока

$$f(x^0) > f(y^0) + \langle \nabla f(y^0), x^0 - y^0 \rangle + \frac{L_0^{j_0}}{2} \|x^0 - y^0\|^2,$$

где

$$x^0 := u^0 := \arg\min_{x \in Q} \varphi_0(x), \ (A_0 :=) \alpha_0 := \frac{1}{L_0^{j_0}},$$

выполнять

$$j_0 := j_0 + 1; \ L_0^{j_0} := 2^{j_0} L_0^0.$$

### Адаптивный Метод Подобных Треугольников

1. $L_{k+1}^0 = L_k^{j_k}/2, \ j_{k+1} = 0.$

2. $\begin{cases} \alpha_{k+1} := \dfrac{1 + A_k \tilde{\mu}}{2L_{k+1}^{j_{k+1}}} + \sqrt{\dfrac{1 + A_k \tilde{\mu}}{4\left(L_{k+1}^{j_{k+1}}\right)^2} + \dfrac{A_k \cdot (1 + A_k \tilde{\mu})}{L_{k+1}^{j_{k+1}}}}, \ A_{k+1} := A_k + \alpha_{k+1}; \\ y^{k+1} := \dfrac{\alpha_{k+1} u^k + A_k x^k}{A_{k+1}}, \ u^{k+1} := \arg\min_{x \in Q} \varphi_{k+1}(x), \ x^{k+1} := \dfrac{\alpha_{k+1} u^{k+1} + A_k x^k}{A_{k+1}}. \end{cases}$ (*)

До тех пор пока

$$f(y^{k+1}) + \langle \nabla f(y^{k+1}), x^{k+1} - y^{k+1} \rangle + \frac{L_{k+1}^{j_{k+1}}}{2} \|x^{k+1} - y^{k+1}\|^2 < f(x^{k+1}),$$

выполнять

$$j_{k+1} := j_{k+1} + 1; \ L_{k+1}^{j_{k+1}} = 2^{j_{k+1}} L_{k+1}^0; \ (*).$$

3. Если не выполнен критерий останова, то $k := k + 1$ и **go to** 1.

В качестве критерия останова, например, можно брать условие

$$\left\| x^{k+1} - \arg\min_{x \in Q} \left\{ \langle \nabla f(x^{k+1}), x - x^{k+1} \rangle + \frac{L_{k+1}^{j_{k+1}}}{2} \|x - x^{k+1}\|^2 \right\} \right\| \leq \tilde{\varepsilon}.$$

Здесь и везде в дальнейшем под "До тех пор пока … выполнять …" подразумевается, что после каждого $j_{k+1} := j_{k+1} + 1$ при следующей проверке условия выхода из этого цикла должным образом меняется не только $L_{k+1}^{j_{k+1}}$, но и $x^{k+1}$, $y^{k+1}$, также входящие в это условие.

**Теорема 6.** *Пусть справедливы предположения 1, 2. Тогда АМПТ для задачи (1) сходится согласно оценке*

$$F\left(x^N\right) - \min_{x \in Q} F(x) \leq \min\left\{\frac{8LR^2}{(N+1)^2}, 2LR^2 \exp\left(-\frac{N}{2}\sqrt{\frac{\mu}{2L\tilde{\omega}_n}}\right)\right\}. \quad (26)$$

*При этом среднее число вычислений значения функции на одной итерации будет $\approx 4$, а градиента функции $\approx 2$.*

**Доказательство.** Нетривиальным в виду оценки (21) и свойства, что все $L_k^{j_k} \leq 2L$, представляется только последняя часть формулировки теоремы. Докажем именно её. Оценим общее число обращений за значениями функции (аналогично получается оценка общего числа обращений за значением градиента функции)

$$\sum_{k=1}^{N} 2(j_k+1) = \sum_{k=1}^{N} 2\left[(j_k-1)+2\right] = \sum_{k=1}^{N} 2\left[\log_2\left(\frac{L_k^{j_k}}{L_{k-1}^{j_{k-1}}}\right)+2\right] = 4N + \log_2\left(\frac{L_N^{j_N}}{L_0^0}\right) \leq 4N + \log_2\left(\frac{2L}{L_0^0}\right).$$

Деля обе части на $N$, получим в правой части приблизительно 4. ∎

В действительности, оценка (26) оказывается, как правило, сильно завышенной, поскольку метод адаптивно настраивается на константу Липшица градиента $L$ на данном участке своего пребывания, а константа $L$, входящая в оценку (26), соответствует (согласно предположению 1) самому плохому случаю (самому плохому участку).

Предположим теперь, что по каким-то причинам невозможно получить точные значения функции и градиента. Тогда соотношение (аналогичное неравенство выписывается и при $k=0$, см. начало этого пункта)

$$f\left(y^{k+1}\right) + \left\langle \nabla f\left(y^{k+1}\right), x^{k+1} - y^{k+1}\right\rangle + \frac{L_{k+1}^{j_{k+1}}}{2}\left\|x^{k+1} - y^{k+1}\right\|^2 \geq f\left(x^{k+1}\right)$$

может не выполниться не при каком $L_{k+1}^{j_{k+1}}$. Допустим, однако, что при этом имеет место

$$f\left(y^{k+1}\right) + \left\langle \nabla f\left(y^{k+1}\right), x^{k+1} - y^{k+1}\right\rangle + \frac{L}{2}\left\|x^{k+1} - y^{k+1}\right\|^2 + \frac{\alpha_{k+1}}{A_{k+1}}\varepsilon \geq f\left(x^{k+1}\right).$$

Тогда заменим в АМПТ соответствующую часть шага 2 на

$$f\left(y^{k+1}\right) + \left\langle \nabla f\left(y^{k+1}\right), x^{k+1} - y^{k+1}\right\rangle + \frac{L_{k+1}^{j_{k+1}}}{2}\left\|x^{k+1} - y^{k+1}\right\|^2 + \frac{\alpha_{k+1}}{A_{k+1}}\varepsilon \geq f\left(x^{k+1}\right). \quad (27)$$

**Теорема 7.** *Пусть справедливо предположение 2 и существует такое число $L > 0$, что любого $k = 1, \ldots, N$ справедливо неравенство*

$$f\left(y^{k+1}\right) + \left\langle \nabla f\left(y^{k+1}\right), x^{k+1} - y^{k+1}\right\rangle + \frac{L}{2}\left\|x^{k+1} - y^{k+1}\right\|^2 + \frac{\alpha_{k+1}}{A_{k+1}}\varepsilon \geq f\left(x^{k+1}\right). \quad (28)$$

*Тогда АМПТ с (27) для задачи (1) сходится согласно оценке*

$$F\left(x^N\right) - \min_{x \in Q} F(x) \leq \min\left\{\frac{8LR^2}{(N+1)^2}, 2LR^2 \exp\left(-\frac{N}{2}\sqrt{\frac{\mu}{2L\tilde{\omega}_n}}\right)\right\} + \varepsilon. \qquad (29)$$

*При этом среднее число вычислений значения функции на одной итерации будет $\approx 4$, а градиента функции $\approx 2$*

**Доказательство.** Ключевым элементом в доказательстве является следующее уточнение леммы 4

$$A_k F\left(x^k\right) \leq \varphi_k^* + A_k \varepsilon, \qquad (30)$$

из которого будет следовать формула (29). Чтобы доказать (30) будем рассуждать по индукции. База индукции $k = 0$ очевидна. Итак, по предположению индукции

$$A_k F\left(x^k\right) - A_k \varepsilon \leq \varphi_k^* = \varphi_k\left(u^k\right),$$

поэтому

$$\varphi_k\left(u^{k+1}\right) \geq \varphi_k\left(u^k\right) + \frac{1 + A_k \tilde{\mu}}{2}\left\|u^{k+1} - u^k\right\|^2 \geq A_k F\left(x^k\right) - A_k \varepsilon + \frac{1 + A_k \tilde{\mu}}{2}\left\|u^{k+1} - u^k\right\|^2.$$

Отсюда

$$\varphi_{k+1}^* \geq A_{k+1} \cdot \left(\frac{A_k}{A_{k+1}} h\left(x^k\right) + \frac{\alpha_{k+1}}{A_{k+1}} h\left(u^{k+1}\right)\right) + A_{k+1} f\left(y^k\right) - A_k \varepsilon$$

$$+ \left\langle \nabla f\left(y^{k+1}\right), \alpha_{k+1} \cdot \left(u^{k+1} - y^{k+1}\right) + A_k \cdot \left(x^k - y^{k+1}\right) \right\rangle + \frac{(1 + A_k \tilde{\mu})}{2}\left\|u^{k+1} - u^k\right\|^2.$$

Следовательно,

$$\varphi_{k+1}^* + A_k \varepsilon \geq A_{k+1}\left[f\left(y^{k+1}\right) + \left\langle \nabla f\left(y^{k+1}\right), x^{k+1} - y^{k+1}\right\rangle + \frac{A_{k+1} \cdot (1 + A_k \tilde{\mu})}{2\alpha_{k+1}^2}\left\|x^{k+1} - y^{k+1}\right\|^2 + h\left(x^{k+1}\right)\right].$$

Отсюда и из условия (27), $A_{k+1} \cdot (1 + A_k \tilde{\mu})/\alpha_{k+1}^2 = L_{k+1}^{j_{k+1}}$ (с учетом (28)) получаем

$$\varphi_{k+1}^* + A_{k+1}\varepsilon = \varphi_{k+1}^* + A_k \varepsilon + \alpha_{k+1}\varepsilon \geq A_{k+1} F\left(x^{k+1}\right), \; L_{k+1}^{j_{k+1}} \leq 2L. \blacksquare$$

В действительности, выше установлено более сильное утверждение – в оценке (29) можно улучшить константу $L$.

**Теорема 8.** *Пусть справедливо предположение 2. Тогда АМПТ с (27) для задачи (1) сходится согласно оценке*

$$F(x^N) - \min_{x \in Q} F(x) \leq \frac{R^2}{A_N} + \varepsilon \leq \min\left\{ \frac{4LR^2}{(N+1)^2}, LR^2 \exp\left(-\frac{N}{2}\sqrt{\frac{\mu}{L\tilde{\omega}_n}}\right) \right\} + \varepsilon, \qquad (31)$$

*где* $L = \max_{k=0,\ldots,N} L_k^{j_k}$. *При этом среднее число вычислений значения функции на одной итерации будет* $\approx 4$, *а градиента функции* $\approx 2$.

Попробуем (подобно [6, 10]) сыграть на условии (27), искусственно вводя неточность.

**Лемма 5** (см. [4, 6]). *Пусть*
$$\|\nabla f(y) - \nabla f(x)\|_* \leq L_\nu \|y - x\|^\nu \qquad (32)$$
*при некотором* $\nu \in [0,1]$. *Тогда*
$$f(y) + \langle \nabla f(y), x - y \rangle + \frac{L}{2}\|x - y\|^2 + \delta \geq f(x), \quad L = L_\nu \left[\frac{L_\nu}{2\delta}\frac{1-\nu}{1+\nu}\right]^{\frac{1-\nu}{1+\nu}}. \qquad (33)$$

Основным результатом данного пункта является описание и последующая оценка скорости сходимости нового варианта универсального метода [6] на базе МПТ (УМПТ).

Положим

$$A_0 = \alpha_0 = 1/L_0^0, \ k = 0, \ j_0 = 0; \ x^0 = u^0 = \arg\min_{x \in Q} \varphi_0(x).$$

До тех пор пока

$$f(x^0) > f(y^0) + \langle \nabla f(y^0), x^0 - y^0 \rangle + \frac{L_0^{j_0}}{2}\|x^0 - y^0\|^2 + \frac{\alpha_0}{2A_0}\varepsilon,$$

где

$$x^0 := u^0 := \arg\min_{x \in Q} \varphi_0(x), \ (A_0 :=)\alpha_0 := \frac{1}{L_0^{j_0}},$$

выполнять

$$j_0 := j_0 + 1; \ L_0^{j_0} := 2^{j_0} L_0^0.$$

### **Универсальный Метод Подобных Треугольников**

> 1. $L_{k+1}^0 = L_k^{j_k}/2, \ j_{k+1} = 0$.

$$2. \begin{cases} \alpha_{k+1} := \dfrac{1+A_k\tilde{\mu}}{2L_{k+1}^{j_{k+1}}} + \sqrt{\dfrac{1+A_k\tilde{\mu}}{4\left(L_{k+1}^{j_{k+1}}\right)^2} + \dfrac{A_k \cdot (1+A_k\tilde{\mu})}{L_{k+1}^{j_{k+1}}}}, \; A_{k+1} := A_k + \alpha_{k+1}; \\ y^{k+1} := \dfrac{\alpha_{k+1}u^k + A_k x^k}{A_{k+1}}, \; u^{k+1} := \arg\min_{x \in Q} \varphi_{k+1}(x), \; x^{k+1} := \dfrac{\alpha_{k+1}u^{k+1} + A_k x^k}{A_{k+1}}. \end{cases} \quad (*)$$

До тех пор пока

$$f(y^{k+1}) + \langle \nabla f(y^{k+1}), x^{k+1} - y^{k+1} \rangle + \frac{L_{k+1}^{j_{k+1}}}{2}\|x^{k+1} - y^{k+1}\|^2 + \frac{\alpha_{k+1}}{2A_{k+1}}\varepsilon < f(x^{k+1}),$$

выполнять

$$j_{k+1} := j_{k+1} + 1; \; L_{k+1}^{j_{k+1}} = 2^{j_{k+1}} L_{k+1}^0; \; (*).$$

3. Если не выполнен критерий останова, то $k := k+1$ и **go to** 1.

**Теорема 9.** *Пусть выполняется условие (32) хотя бы для $\nu = 0$, и справедливо предположение 2 с $\mu \geq 0$ (допускается брать $\mu = 0$). Тогда УМПТ для задачи (1) сходится согласно оценке*

$$F(x^N) - \min_{x \in Q} F(x) \leq \varepsilon,$$

$$N \approx \min\left\{ \inf_{\nu \in [0,1]} \left( \frac{L_\nu \cdot (16R)^{1+\nu}}{\varepsilon} \right)^{\frac{2}{1+3\nu}}, \inf_{\nu \in [0,1]} \left\{ \left( \frac{8L_\nu^{\frac{2}{1+\nu}} \tilde{\omega}_n}{\mu \varepsilon^{\frac{1-\nu}{1+\nu}}} \right)^{\frac{1+\nu}{1+3\nu}} \ln^{\frac{2+2\nu}{1+3\nu}} \left( \frac{16 L_\nu^{\frac{4+6\nu}{1+\nu}} R^2}{(\mu/\tilde{\omega}_n)^{\frac{1+\nu}{1+3\nu}} \varepsilon^{\frac{5+7\nu}{2+6\nu}}} \right) \right\} \right\}. \quad (34)$$

*При этом среднее число вычислений значения функции на одной итерации будет $\approx 4$, а градиента функции $\approx 2$.*

**Доказательство.** Рассмотрим два случая, когда $\mu \geq 0$ – мало: $\mu \ll \varepsilon/(2R^2)$, $\mu$ – велико: $\mu \gg \varepsilon/(2R^2)$, см. формулу (23).

В первом случае будем считать, что

$$A_{k+1}/\alpha_{k+1}^2 \approx A_{k+1} \cdot (1+A_k\tilde{\mu})/\alpha_{k+1}^2 = L_{k+1}^{j_{k+1}}, \text{ т.е. } \varepsilon \frac{\alpha_{k+1}}{2A_{k+1}} \approx \frac{\varepsilon}{2}\sqrt{\frac{1}{L_{k+1}^{j_{k+1}} A_{k+1}}}, \quad (35)$$

а во втором случае

$$A_{k+1}^2 \tilde{\mu}/\alpha_{k+1}^2 \approx A_{k+1} \cdot (1+A_k\tilde{\mu})/\alpha_{k+1}^2 = L_{k+1}^{j_{k+1}}, \text{ т.е. } \varepsilon \frac{\alpha_{k+1}}{2A_{k+1}} \approx \frac{\varepsilon}{2}\sqrt{\frac{\tilde{\mu}}{L_{k+1}^{j_{k+1}}}}. \quad (36)$$

Из формулы (31) (см. теорему 8) имеем

$$\frac{R^2}{A_N} + \frac{\varepsilon}{2} \approx \varepsilon,$$

т.е. $A_N \approx 2R^2/\varepsilon$, а также ( $L = \max_{k=0,\ldots,N} L_k^{j_k}$ )

$$N^2 \approx \frac{8LR^2}{\varepsilon}, \quad \text{(в первом случае)} \tag{37}$$

$$N^2 \approx 4\frac{L}{\tilde{\mu}}\ln^2\left(\frac{2LR^2}{\varepsilon}\right). \quad \text{(во втором случае)} \tag{38}$$

Из формул (33), (35) – (38) имеем, что в первом случае

$$L \le 2L_\nu \left[\frac{L_\nu}{2\frac{\varepsilon}{2}\sqrt{\frac{1}{LA_N}}}\frac{1-\nu}{1+\nu}\right]^{\frac{1-\nu}{1+\nu}} \le 2L_\nu\left[\frac{L_\nu\sqrt{LA_N}}{\varepsilon}\frac{1-\nu}{1+\nu}\right]^{\frac{1-\nu}{1+\nu}} \le 2L_\nu^{\frac{2}{1+\nu}}\left[\frac{N}{2\varepsilon}\frac{1-\nu}{1+\nu}\right]^{\frac{1-\nu}{1+\nu}}, \tag{39}$$

а во втором случае

$$L \le 2L_\nu\left[\frac{L_\nu}{2\frac{\varepsilon}{2}\sqrt{\frac{\tilde{\mu}}{L}}}\frac{1-\nu}{1+\nu}\right]^{\frac{1-\nu}{1+\nu}} \le 2L_\nu\left[\frac{L_\nu\sqrt{L/\tilde{\mu}}}{\varepsilon}\frac{1-\nu}{1+\nu}\right]^{\frac{1-\nu}{1+\nu}} \le 2L_\nu^{\frac{2}{1+\nu}}\left[\frac{N}{2\varepsilon}\frac{1-\nu}{1+\nu}\right]^{\frac{1-\nu}{1+\nu}}. \tag{40}$$

Подставляя (39) в (37), а (40) в (38) и учитывая, что параметр $\nu \in [0,1]$ можно выбирать произвольно (допускается, что $L_\nu = \infty$ при некоторых $\nu$ – важно, чтобы существовало хотя бы одно значение $\nu$ при котором $L_\nu < \infty$; по условию $L_0 < \infty$), получим соответственно,

$$N^2 \approx \frac{16L_\nu^{\frac{2}{1+\nu}}\left[\frac{N}{2\varepsilon}\frac{1-\nu}{1+\nu}\right]^{\frac{1-\nu}{1+\nu}}R^2}{\varepsilon} \Rightarrow N^{\frac{1+3\nu}{1+\nu}} \approx \frac{16L_\nu^{\frac{2}{1+\nu}}R^2}{\varepsilon^{\frac{2}{1+\nu}}} \Rightarrow N \approx \inf_{\nu \in [0,1]}\left(\frac{L_\nu \cdot (16R)^{1+\nu}}{\varepsilon}\right)^{\frac{2}{1+3\nu}}, \tag{41}$$

$$N^2 \approx \frac{8L_\nu^{\frac{2}{1+\nu}}\left[\frac{N}{2\varepsilon}\frac{1-\nu}{1+\nu}\right]^{\frac{1-\nu}{1+\nu}}}{\tilde{\mu}}\ln^2\left(\frac{2L_\nu^2 R^2 N}{\varepsilon^{3/2}}\right) \Rightarrow N^{\frac{1+3\nu}{1+\nu}} \approx \frac{8L_\nu^{\frac{2}{1+\nu}}}{\tilde{\mu}\varepsilon^{\frac{1-\nu}{1+\nu}}}\ln^2\left(\frac{2L_\nu^2 R^2 N}{\varepsilon^{3/2}}\right) \Rightarrow$$

$$\Rightarrow N \approx \left( \frac{8L_\nu^{\frac{2}{1+\nu}}}{\tilde{\mu}\varepsilon^{\frac{1-\nu}{1+\nu}}} \right)^{\frac{1+\nu}{1+3\nu}} \ln^{\frac{2+2\nu}{1+3\nu}} \left( \frac{2L_\nu^2 R^2 N}{\varepsilon^{3/2}} \right) \Rightarrow N \approx \inf_{\nu \in [0,1]} \left\{ \left( \frac{8L_\nu^{\frac{2}{1+\nu}}}{\tilde{\mu}\varepsilon^{\frac{1-\nu}{1+\nu}}} \right)^{\frac{1+\nu}{1+3\nu}} \ln^{\frac{2+2\nu}{1+3\nu}} \left( \frac{16L_\nu^{\frac{4+6\nu}{1+\nu}} R^2}{\tilde{\mu}^{\frac{1+\nu}{1+3\nu}} \varepsilon^{\frac{5+7\nu}{2+6\nu}}} \right) \right\}. \quad (42)$$

Из формул (41), (42) получаем оценку (34). ∎

Более аккуратные рассуждения позволяют в несколько раз уменьшить константы, входящие в оценку (34). Оценка (34) согласуется с нижними оценками для соответствующих классов задач [20].

Практические эксперименты, проведенные (А.Ю. Горновым) с УМПТ (при $\mu = 0$) показали, что этот метод работает подобно другим версиям универсального метода [6, 10]. К сожалению, и в такой (предложенной в статье – УМПТ) модификации метод типично проигрывает на практике методу сопряженных градиентов на квадратичных задачах безусловной оптимизации. Впрочем, то что обычные (не универсальные) быстрые градиентные методы проигрывают на таких задачах методу сопряженных градиентов было известно и ранее [21]. Новое здесь то, что универсальные быстрые градиентные методы также проигрывают.

**Замечание 4 (А.И. Тюрин).** Все описанные методы могут быть распространены (насколько нам известно, это пока еще не сделано в общем случае) на задачи условной оптимизации. Для этого сначала, следуя п. 2.3 [18], стоит рассмотреть минимаксную задачу (ввести правильную лианеризацию исходного функционала и градиентное отображение). Далее использовать идею метода нагруженных функционалов п. 2.3.4 [18], приводящую к рестартам по неизвестному параметру п. 2.3.5 [18] (оптимальное значение функционала задачи), введение которого, позволяет свести задачу условной оптимизации к минимаксной. Дополнительная плата за такое "введение" (т.е. за рестарты) будет всего лишь логарифмическая, и с точностью до этой "платы" оценки (на число итераций) будут оптимальными. Однако, к сожалению, итерации метода получаются слишком дорогими. В отличие от п. 2.3 [18] сложность итерации растет с ростом её номера. От этого можно избавиться (с точностью до рестартов), используя вместо линейки описанных в этой статье методов, например, быстрый градиентный метод в форме Allen-Zhu–Orecchia [10]. Другой способ, заменить в Методе Подобных Треугольников $u^{k+1} = \arg\min_{x \in Q} \varphi_{k+1}(x)$ на

$$u^{k+1} = \arg\min_{x \in Q} \left\{ \alpha_{k+1} \cdot \left( \langle \nabla f(y^{k+1}), x - y^{k+1} \rangle + h(x) \right) + V(x, u^k) \right\}.$$

Детали см. в работах [22, 23].

**Замечание 5 (А.И. Тюрин).** Все описанные методы могут быть распространены на случай, когда шаг $u^{k+1} = \arg\min_{x \in Q} \varphi_{k+1}(x)$ может быть осуществлен с погрешностью. Соответствующие выкладки практически дословно повторяют рассуждения п. 5.5 [11]. Детали см. в работе [22].

**Замечание 6.** В большинстве приложений "стоимость" (время) получения от оракула (роль которого, как правило, играют нами же написанные подпрограммы вычисления градиента) градиента функционала заметно превышает время, затрачиваемое на то, чтобы сделать шаг итерации, исходя из выданного оракулом вектора. Желание сбалансировать это рассогласование (усложнить итерации, сохранив при этом старый порядок сложности, и выиграть за счет этого в сокращении числа итераций), привело к возникновению композитной оптимизации [2], в которой (аддитивная) часть функционала задачи переносится без лианеризации (запроса градиента) в итерации. Другой способ перенесения части сложности задачи на итерации был описан в замечании 4. Здесь остается еще много степеней свободы, позволяющих играть на том насколько "дорогой" будет оракул и соответствующая (этому оракулу) "процедура проектирования", и том сколько (внешних) итераций потребуется методу для достижения заданной точности. В частности, если обращение к оракулу за градиентом и последующее проектирование требуют, в свою очередь, решения вспомогательных оптимизационных задач, то можно "сыграть" на том, насколько точно надо решать эти вспомогательные задачи, пытаясь найти "золотую середину" между стоимостью итерации и числом итераций (см. замечание 5). Также можно сыграть и на том, как выделять эти вспомогательные задачи. Другими словами, что понимать под оракулом и под итерацией метода. Общая идея "разделяй и властвуй", применительно к численным методам выпуклой оптимизации может принимать довольно неожиданные и при этом весьма эффективные формы (ярким примером являются методы внутренней точки [11, 18]). Разные варианты описанной игры в связи с транспортно-сетевыми приложениями уже разбирались нами в других работах [9, 24, 25].

## 5. Универсальный метод подобных треугольников для задач стохастической композитной оптимизации

Предположим теперь, что вместо настоящих градиентов нам доступны только стохастические градиенты $\nabla f(x) \to \nabla f(x,\xi)$ (см., например, [3, 8]).

Для большей наглядности далее в этом разделе удобно будет считать, что $\|\ \| = \|\ \|_2$.

**Предположение 3.** *Для всех* $x \in Q$

$$E_\xi\left[\nabla f(x,\xi)\right] = \nabla f(x) \ \text{и} \ E_\xi\left[\left\|\nabla f(x,\xi) - \nabla f(x)\right\|_*^2\right] \le D. \tag{43}$$

В ряде приложений бывает полезно рассматривать модификацию условия (43)

$$LE_\xi\left[\max_{x,y \in Q}\left\{\langle \nabla f(y,\xi) - \nabla f(y), x-y\rangle - \frac{L}{2}\|x-y\|^2\right\}\right] \le \tilde{D}.$$

Далее приводится стохастический вариант УМПТ (СУМПТ). По-видимому, это первая попытка перенесения универсального метода на задачи стохастической оптимизации.

Предварительно введём обозначения

$$\overline{\nabla}^m f(x) = \frac{1}{m}\sum_{k=1}^{m}\nabla f(x,\xi^k), \qquad (44)$$

где $\xi^k$ – независимые одинаково распределенные (так же как $\xi$) случайные величины. В принципе, можно было бы аналогично ввести

$$\overline{f}^{\tilde{m}}(x) = \frac{1}{\tilde{m}}\sum_{k=1}^{\tilde{m}} f(x,\xi^k)$$

и распространить приведенные далее результаты на случай, когда и значение функции $f(x)$ необходимо оценивать. Однако для большей наглядности ограничимся далее случаем, когда значение $f(x)$ можно точно посчитать (т.е. нет необходимости его оценивать).

Переопределим последовательность (14)

$$\varphi_0(x) = V(x,y^0) + \alpha_0\left[f(y^0) + \langle \overline{\nabla}^m f(y^0), x - y^0\rangle + \tilde{\mu}V(x,y^0) + h(x)\right],$$

$$\varphi_{k+1}(x) = \varphi_k(x) + \alpha_{k+1}\left[f(y^{k+1}) + \langle \overline{\nabla}^m f(y^{k+1}), x - y^{k+1}\rangle + \tilde{\mu}V(x,y^k) + h(x)\right].$$

Положим

$$A_0 = \alpha_0 = 1/L_0^0,\ m_0 := \frac{2DA_0}{L_0^0 \alpha_0 \varepsilon},\ k = 0,\ j_0 = 0;\ x^0 = u^0 = \arg\min_{x\in Q}\varphi_0(x).$$

До тех пор пока

$$f(x^0) > f(y^0) + \langle \overline{\nabla}^{m_0} f(y^0), x^0 - y^0\rangle + \frac{L_0^{j_0}}{2}\|x^0 - y^0\|^2 + \frac{3\alpha_{k+1}}{2A_{k+1}}\varepsilon,$$

где

$$x^0 := u^0 := \arg\min_{x\in Q}\varphi_0(x) =$$

$$= \arg\min_{x\in Q}\left\{V(x,y^0) + \alpha_0\left[f(y^0) + \langle \overline{\nabla}^{m_0} f(y^0), x - y^0\rangle + \tilde{\mu}V(x,y^0) + h(x)\right]\right\},$$

$$(A_0 :=)\alpha_0 := \frac{1}{L_0^{j_0}},\ m_0 := \frac{2DA_0}{L_0^{j_0}\alpha_0 \varepsilon},$$

выполнять

$$j_0 := j_0 + 1;\ L_0^{j_0} := 2^{j_0}L_0^0.$$

### **Стохастический Универсальный Метод Подобных Треугольников**

1. $L_{k+1}^0 = L_k^{j_k}/2,\ j_{k+1} = 0.$

2. $\begin{cases} \alpha_{k+1} := \dfrac{1+A_k\tilde{\mu}}{2L_{k+1}^{j_{k+1}}} + \sqrt{\dfrac{1+A_k\tilde{\mu}}{4\left(L_{k+1}^{j_{k+1}}\right)^2} + \dfrac{A_k\cdot(1+A_k\tilde{\mu})}{L_{k+1}^{j_{k+1}}}},\ A_{k+1} := A_k + \alpha_{k+1},\ m_{k+1} := \dfrac{2DA_{k+1}}{L_{k+1}^{j_{k+1}}\alpha_{k+1}\varepsilon}; \\ y^{k+1} := \dfrac{\alpha_{k+1}u^k + A_k x^k}{A_{k+1}},\ u^{k+1} := \arg\min_{x\in Q}\varphi_{k+1}(x),\ x^{k+1} := \dfrac{\alpha_{k+1}u^{k+1} + A_k x^k}{A_{k+1}}. \end{cases}$ (*)

До тех пор пока

$$f(y^{k+1}) + \langle \overline{\nabla}^{m_{k+1}} f(y^{k+1}), x^{k+1} - y^{k+1}\rangle + \frac{L_{k+1}^{j_{k+1}}}{2}\|x^{k+1} - y^{k+1}\|^2 + \frac{3\alpha_{k+1}}{2A_{k+1}}\varepsilon < f(x^{k+1}),$$

выполнять

$$j_{k+1} := j_{k+1} + 1;\ L_{k+1}^{j_{k+1}} = 2^{j_{k+1}} L_{k+1}^0;\ (*).$$

3. Если не выполнен критерий останова, то $k := k+1$ и **go to** 1.

**Теорема 10.** *Пусть выполняется условие (32) хотя бы для $\nu = 0$, справедливо предположение 2 с $\mu \geq 0$ (допускается брать $\mu = 0$), справедливо предположение 3. Тогда СУМПТ для задачи (1) сходится согласно оценке*

$$E\left[F(x^N)\right] - \min_{x\in Q} F(x) \leq 2\varepsilon,$$

$$N \approx \min\left\{\inf_{\nu\in[0,1]}\left(\frac{L_\nu\cdot(32R)^{1+\nu}}{\varepsilon}\right)^{\frac{2}{1+3\nu}},\ \inf_{\nu\in[0,1]}\left\{\left(\frac{16L_\nu^{\frac{2}{1+\nu}}\tilde{\omega}_n}{\mu\varepsilon^{\frac{1-\nu}{1+\nu}}}\right)^{\frac{1+\nu}{1+3\nu}}\ln^{\frac{2+2\nu}{1+3\nu}}\left(\frac{32L_\nu^{\frac{4+6\nu}{1+\nu}}R^2}{(\mu/\tilde{\omega}_n)^{\frac{1+\nu}{1+3\nu}}\varepsilon^{\frac{5+7\nu}{2+6\nu}}}\right)\right\}\right\}. \quad (45)$$

*Оценка (45) – это оценка числа итераций. При этом среднее число вычислений значения функции на одной итерации будет $\approx 4$. Однако не менее интересна оценка числа обращений за стохастическим градиентом*

$$Q \approx 2\min\left\{\frac{4DR^2}{\varepsilon^2},\ \frac{2D\tilde{\omega}_n}{\mu\varepsilon}\ln\left(\frac{2L_0^{j_0}R^2}{\varepsilon}\right)\right\} + 2N. \quad (46)$$

**Доказательство.** По неравенству Фенхеля (см. главу 7 [3])

$$\langle \overline{\nabla}^{m_{k+1}} f(y^{k+1}) - \nabla f(y^{k+1}), x^{k+1} - y^{k+1}\rangle - \frac{L_{k+1}^{j_{k+1}}/2}{2}\|x^{k+1} - y^{k+1}\|^2 \leq \frac{2}{L_{k+1}^{j_{k+1}}}\|\overline{\nabla}^{m_{k+1}} f(y^{k+1}) - \nabla f(y^{k+1})\|_*^2.$$

Поэтому ключевое неравенство (27)

$$f(y^{k+1}) + \langle \overline{\nabla}^{m_{k+1}} f(y^{k+1}), x^{k+1} - y^{k+1}\rangle + \frac{L_{k+1}^{j_{k+1}}}{2}\|x^{k+1} - y^{k+1}\|^2 + \frac{\alpha_{k+1}}{A_{k+1}}\frac{\varepsilon}{2} \geq f(x^{k+1})$$

переписывается следующим образом

$$f\left(y^{k+1}\right)+\left\langle \nabla f\left(y^{k+1}\right), x^{k+1}-y^{k+1}\right\rangle + L_{k+1}^{j_{k+1}}\left\|x^{k+1}-y^{k+1}\right\|^2 + \frac{\alpha_{k+1}}{A_{k+1}}\frac{\varepsilon}{2} + \frac{\alpha_{k+1}}{A_{k+1}}\varepsilon \geq f\left(x^{k+1}\right). \quad (47)$$

При получении неравенства (47) для большей наглядности заменяем в рассуждениях правую часть в неравенстве Фенхеля оценкой его математического ожидания, равной согласно условию (43)

$$\frac{2D}{L_{k+1}^{j_{k+1}}m_{k+1}} = \frac{\alpha_{k+1}}{A_{k+1}}\varepsilon.$$

В действительности, тут требуются более громоздкие рассуждения, см., например, [3, 8, 22], которые приведут к необходимости увеличения в несколько раз (во сколько именно раз, зависит от "тяжести" хвостов распределения $\nabla f(x,\xi)$ и от выбранного доверительного уровня) константы 2 в формуле выбора $m_{k+1}$ в описании СУМПТ, и к соответствующему увеличению $N$ и Q.

Оценим число обращений за стохастическим градиентом Q, используя схему доказательства теоремы 9. Для этого, прежде всего, заметим, что $A_N \approx 2R^2/\varepsilon$. Рассмотрим два случая, когда $\mu \geq 0$ – мало: $\mu \ll \varepsilon/(2R^2)$, $\mu$ – велико: $\mu \gg \varepsilon/(2R^2)$.

В первом случае будем считать, что (см. формулу (35))

$$A_{k+1}/\alpha_{k+1}^2 \approx A_{k+1}\cdot\left(1+A_k\tilde{\mu}\right)/\alpha_{k+1}^2 = L_{k+1}^{j_{k+1}}.$$

а во втором случае, что (см. формулу (36))

$$A_{k+1}^2\tilde{\mu}/\alpha_{k+1}^2 \approx A_{k+1}\cdot\left(1+A_k\tilde{\mu}\right)/\alpha_{k+1}^2 = L_{k+1}^{j_{k+1}}.$$

В первом случае число обращений за стохастическим градиентом оценивается, соответственно, как

$$Q \approx \sum_{k=0}^{N}\frac{2DA_k}{L_k^{j_k}\alpha_k\varepsilon} = \frac{2D}{\varepsilon}\sum_{k=0}^{N}\frac{A_k}{L_k^{j_k}\alpha_k} \approx \frac{2D}{\varepsilon}\sum_{k=0}^{N}\alpha_k = \frac{2D}{\varepsilon}A_N \approx \frac{4DR^2}{\varepsilon^2}, \quad (48)$$

$$Q \approx \sum_{k=0}^{N}\frac{2DA_k}{L_k^{j_k}\alpha_k\varepsilon} = \frac{2D}{\varepsilon}\sum_{k=0}^{N}\frac{A_k}{L_k^{j_k}\alpha_k} \approx \frac{2D}{\tilde{\mu}\varepsilon}\sum_{k=0}^{N}\frac{\alpha_k}{A_k} \approx \frac{2D}{\tilde{\mu}\varepsilon}\int_{1}^{A_N/\alpha_0}\frac{dA}{A} \approx \frac{2D}{\tilde{\mu}\varepsilon}\ln\left(\frac{A_N}{\alpha_0}\right) \approx \frac{2D}{\tilde{\mu}\varepsilon}\ln\left(\frac{2L_0^{j_0}R^2}{\varepsilon}\right). \quad (49)$$

Из формул (48), (49) получаем оценку (46). ∎

Оценка (46) с точностью до логарифмических множителей соответствует нижней оценки [12]. В случае $\|\ \|\neq\|\ \|_2$ в приведенных выше оценках может возникать дополнительный множитель $\sim \ln n$, $n = \dim x$ (Proposition 6, [26]).

**Замечание 7.** Все приведенные в статье результаты допускают обобщение на случай, когда вместо точного оракула, выдающего значения функции $f(x)$ и ее градиента $\nabla f(x)$ (или их несмещенные реализации), используется (стохастический) $(\delta, L, \mu)$-оракул с $\delta = \mathrm{O}(\varepsilon/N)$ и $L = \mathrm{O}\left(\max_{k=0,\ldots,N} L_k^{j_k}\right)$ [3, 4, 7–10]. Такое обобщение существенным образом используется, например, в замечании 6 выше.